\documentclass[12pt]{amsart}
\usepackage{amssymb}
\theoremstyle{plain}
\newtheorem
{thm}{Theorem}[section]
\newtheorem
{proposition}[thm]{Proposition}
\newtheorem
{lemma}[thm]{Lemma}

\newtheorem
{corollary}[thm]{Corollary}
\newcommand{\be}{\begin{equation*}}
\newcommand{\ee}{\end{equation*}}
\newcommand{\beq}{\begin{equation}}
\newcommand{\eeq}{\end{equation}}
\newcommand{\bg}{\begin{gather*}}
\newcommand{\eg}{\end{gather*}}
\newcommand{\bga}{\begin{gather}}
\newcommand{\ega}{\end{gather}}

\newcommand{\ka}{\kappa}
\newcommand{\Gm}{\Gamma}
\newcommand{\la}{\lambda}
\newcommand{\thi}{\vartheta}
\newcommand{\Ref}[1]{{(\ref{#1})}}
\newcommand{\tr}{\operatorname{Tr}}
\newcommand{\End}{\operatorname{End}}

\renewcommand{\Im}{\operatorname{Im}}
\theoremstyle{definition}

\newtheorem*{example}{Example}
\newtheorem*{remark}{Remark}
\title[Modular transformations]
{Modular transformations of the elliptic hypergeometric functions,
Macdonald polynomials, and the shift operator}
\author[G. Felder, L. Stevens, and A. Varchenko]
{G. Felder $^{\,\star}$,
L. Stevens$^{\, \diamond}$, and
A. Varchenko$^{\,\diamond,1}$}
\thanks{$^1$ Supported in part by NSF grant DMS-9801582}
\begin{document}
\begin{abstract}
We consider the space of elliptic hypergeometric functions of the
$sl_{2}$ type associated with elliptic curves with one marked point. This
space represents conformal blocks in the $sl_2$ WZW model of CFT. The
modular group acts on this space. We give formulas for the matrices of
the action in terms of values at roots of unity of Macdonald polynomials
of the $sl_2$ type.
\end{abstract}
\maketitle
\begin{center}
{\it
$^\star$ Departement Mathematik, ETH-Zentrum, 8092 Z\"urich, Switzerland,

felder@math.ethz.ch

\medskip

$^\diamond$ Department of Mathematics, University of North Carolina
at Chapel Hill,

Chapel Hill, NC 27599-3250, USA,

stevens@math.unc.edu, av@math.unc.edu}
\end{center}

\centerline{February, 2002}

\bigskip
\centerline{\emph{Dedicated to V.~I.~Arnold on his 65th birthday}}



\section{Introduction} In the WZW model of conformal field theory associated
with a simple complex Lie algebra $\mathfrak{g}$,
one defines a holomorphic
vector bundle, the bundle of conformal blocks,
on the moduli space of smooth complex compact curves with
marked points labeled by representations of $\mathfrak{g}$.  This vector bundle comes with a projectively flat connection,
see \cite{TUY}. For curves of genus zero and one, this
 connection is flat and can be described in explicit classical terms.
Moreover, horizontal sections admit
integral representations
\cite{SV}, \cite{FV1}.

For genus zero curves, the connection
is the Knizhnik--Zamolodchikov (KZ) connection and the
equation for
horizontal sections (KZ equation) is a generalization of the
Gauss hypergeometric equation. The KZ equation reduces to the Gauss
hypergeometric equation in a special case.
The solutions to the KZ equations are given by integrals
of  certain differential forms over
 $m$-dimensional cycles
 where $m$ depends on $\mathfrak g$ and on the representations
  labeling marked points. The classical
  integral representation of the hypergeometric function
  is recovered in the case of four points on the
Riemann sphere when $m=1$.
The horizontal sections (conformal blocks)
are holomorphic functions on the
universal covering of the configuration space of points
on the complex plane and integral representation may be
used to compute the action of covering transformations
on conformal blocks. In this way representations of the
pure braid group are obtained as monodromy representations.

We consider here the case of genus one curves with one
marked point (elliptic curves)
and the Lie algebra $sl_2$. The equation for
horizontal sections is the
Kni\-zhnik--Za\-mo\-lo\-dchi\-kov--Ber\-nard (KZB) heat
equation \eqref{KZB}: it is essentially the heat equation
associated to the Hamilton operator of a particle
in a Weierstrass function potential. Again, solutions are
given by generalizations of hypergeometric integrals,
which are appropriately called elliptic hypergeometric
integrals. Of particular interest is a finite dimensional
subspace of the space of solutions, the space of conformal
blocks of the WZW model. It may be characterized by
symmetry and holomorphy conditions, see Sect.~2. This
subspace is invariant under the (projective)
  action of the modular group
  $\mathrm{SL}(2,\mathbb{Z})$ of covering transformations
  of the upper half plane, viewed as the universal cover of
  the moduli space of elliptic curves.

  We compute the projective
  action of the modular group on this space
  and relate the matrix elements of the $S$-transformation
  $\tau\mapsto-1/\tau$ to values of $A_1$-Macdonald polynomials. This
  implies a
  special case of Kirillov's theorem, stating that the
  representation of the modular group on certain conformal blocks
 on elliptic curves
for $sl_N$ is equivalent to a representation where the
  matrix elements of $S$ are given in terms of
  $A_{N-1}$-Macdonald polynomials.

  One new feature of our computation is that it gives
  a formula for the matrix elements of the
  projective representation of the
  modular group (not just its conjugacy class)
  with respect to an explicit basis of
  elliptic hypergeometric integral solutions of the KZB equation
  \eqref{KZB}.
  Moreover, the elliptic hypergeometric integrals give
  naturally a recursive procedure, developed in \cite{FS}
  in a similar situation,
  to compute the matrix elements of $S$ by repeated application
  of the Stokes theorem.  This recursive
  procedure shows the role of the {\em shift operator} \cite{AI}
  in this context: it is identical to the recursive construction
of Macdonald polynomials out of Schur functions
  by repeated application of the shift operator.
  We also discuss the relation of the matrix elements of $S$ with
  the traces of intertwining operators of the quantum group
  $U_q(sl_2)$ at root of unity.

  Finally, let us point out that the KZB connection is unitary
  with respect to a hermitian form which can also
  be given by integrals of elliptic hypergeometric type.
  They are discussed in \cite{FG}.

The authors thank P. Etingof and A. Kirillov, Jr., for useful discussions
and the referee for useful suggestions.

\section{Conformal Blocks on the Torus}
Let $\ka$ and $p$ be non-negative integers satisfying $\ka\geq2p+2$.   Let
$q=e^{\frac{\pi i}{\ka}}$.  Denote
\begin{gather*}
[n]_{q}=[n]=\frac{q^n-q^{-n}}{q-q^{-1}},\qquad
[n]!=[1][2]\cdots[n], \qquad
\begin{bmatrix}n\\j\end{bmatrix}=\frac{[n]!}{[j]![n-j]!},\\
(n,q)_{j}=[n][n+1]\cdots[n+j-1].
\end{gather*}
Let $\tau\in\mathbb{C}$ be such that $\Im\tau>0$.
The KZB-heat equation is the partial differential equation
\begin{equation}\label{KZB}
2\pi i\ka{\partial u\over\partial\tau}(\la,\tau)={\partial^2
  u\over\partial\la^{2}}(\la,\tau)+p(p+1)\rho\,'(\la,\tau)u(\la,\tau).
\end{equation}
Here, the prime denotes the derivative with respect to the first
argument, and $\rho$ is defined in terms of 
the first Jacobi theta function \cite{WW},
\begin{equation*}
\thi_{1}(t,\tau)=-\sum_{j\in\mathbb{Z}}
e^{\pi i(j+\frac{1}{2})^2\tau+2\pi i(j+\frac{1}{2})(t+\frac{1}{2})},\qquad
\rho(t,\tau)=\frac{\thi_{1}'(t,\tau)}{\thi_{1}(t,\tau)}.
\end{equation*}
Holomorphic solutions of the KZB-heat equation with the properties,
\begin{enumerate}
\renewcommand{\theenumi}{\roman{enumi}}\renewcommand{\labelenumi}{(\theenumi)}
\item $u(\la+2,\tau)=u(\la,\tau)$,
\item $u(\la+2\tau,\tau)=e^{-2\pi i\ka(\la+\tau)}\,u(\la,\tau)$, 
\item $u(-\la,\tau)=(-1)^{p+1}\,u(\la,\tau)$,
\item $u(\la,\tau)=\mathcal{O}((\la-m-n\tau)^{p+1})$ as $\la\to
  m+n\tau$ for any $m,n\in\mathbb{Z}$,
\end{enumerate}
are called conformal blocks associated with the family of elliptic curves
$\mathbb{C}/\mathbb{Z}+\tau\mathbb{Z}$ with the marked point $z=0$ and the
irreducible $sl_{2}$ representation of dimension $2p+1$.  It is known
that the space of conformal blocks has dimension $\ka-2p-1$.

{\bf Remark.}
The relation of this definition of conformal blocks
with the more standard one of \cite{TUY} is the following. 
Recall the definition of \cite{TUY} in this case:
integrable irreducible representations $V_{k,m}$
of the affine Kac--Moody Lie algebra $\widehat{sl}_2$ 
are labeled by two non-negative integers: $k$, the level and
$m\in\{0,\dots,k\}$, the highest weight. To each such pair $(k,m)$
and a curve $E_\tau=\mathbb{C}/\mathbb{Z}+\tau\mathbb{Z}$ one associates
a space of conformal blocks: it
is the space of linear functions on $V_{k,m}$ invariant under a
natural action
the Lie algebra $R(\tau)$ of $sl_2$-valued functions on $E_\tau$
with poles at $0$. As one varies $\tau$ in the upper half plane,
the spaces of conformal blocks form a vector bundle with a flat
connection. The description of this bundle and connection in
terms of theta functions was explained in \cite{FW}: 
the family of Lie algebras
$R(\tau)$ is extended to a two parameter family $R(\lambda,\tau)$
which also acts on $V_{k,m}$. The $R(\lambda,\tau)$-invariant linear 
forms on $V_{k,m}$ form
a vector bundle of twisted conformal blocks with a flat connection.
The space of horizontal sections is then isomorphic to the space
of twisted conformal blocks at any $(\lambda,\tau)$, in particular
if $\lambda=0$, where (untwisted) conformal blocks are recovered.
We have an injective restriction map from the space of horizontal sections
of the bundle of twisted conformal blocks
to the space of holomorphic functions of $\lambda,\tau$ with
values in the the dual of the zero weight space $V_m[0]\subset V_{k,m}$ of the
$2m+1$-dimensional irreducible representation $V_m$ of $sl_2$. If
$m$ is odd, $V_m[0]=0$. If $m=2p$ is even, then $V_m[0]=\mathbb{C}$ and the
image of the restriction map
consists of functions of the form $v(\lambda,\tau)=
\vartheta_1(\lambda,\tau)u(\lambda,\tau)$, where $u$ is a solution 
of the KZB-heat equation obeying (i)-(iv) with $\kappa=k+2$.
Conformal blocks may also be defined in other ways which also
lead to the same description in terms of theta functions: 
in \cite{FG} this is done for conformal blocks defined
as states of the Chern--Simons--Witten theory and
in \cite{EFK} for the definition
as spherical functions on the Kac--Moody group.

Introduce two transformations
\begin{equation}\label{defnsofst}
Tu(\la,\tau)=u(\la,\tau+1),\qquad 
Su(\la,\tau)=e^{-\pi i \ka{\la^{2}\over2\tau}}\tau^{-\frac{1}{2}-{p(p+1)\over\ka}}\,u\left({\la\over\tau},-{1\over\tau}\right),
\end{equation}
where we fix $\arg\tau\in(0,\pi)$.  
\begin{thm}\label{kzbsols}\cite{EK}
If $u(\la,\tau)$ is a solution of the KZB-heat equation, then
$Tu(\la,\tau)$ and $Su(\la,\tau)$
are solutions too.  Moreover, the transformations $T$ and $S$ preserve
the properties (i)-(iv).  
\end{thm}
The proof of Theorem \ref{kzbsols} is by direct verification.
Restricted to the space of
conformal blocks, the transformations $T$ and $S$ satisfy the relations
\begin{equation*}
S^{2}=(-1)^{p}\,iq^{-p(p+1)}I,\qquad (ST)^{3}=(-1)^{p}\,iq^{-p(p+1)}I,
\end{equation*}  
where $I$ is the identity transformation.

The modular group is the group generated by two elements $T$ and $S$
with relations
\begin{equation*}
S^{2}=1,\qquad (ST)^{3}=1.
\end{equation*}
The modular group is naturally isomorphic to the quotient group
$\mathrm{SL}(2,\mathbb{Z})/\{\pm I\}$.
The formulas \Ref{defnsofst} define a
projective representation of the modular group in the space of
conformal blocks.

The isomorphism class of the projective representation of the modular
group in the space of conformal blocks is described in \cite{K} using
the Kazhdan-Lusztig-Finkelberg isomorphism between the modular tensor
categories arising from affine Lie algebras
and the modular tensor categories arising from quantum groups at roots of unity.
In this paper, we describe the action of the modular group on the
space of conformal blocks using integral representations of solutions
of the KZB-heat equation.
\section{Integral Representations of Solutions of the KZB-Heat Equation} 
Solutions of the KZB-heat equation can be realized as elliptic hypergeometric
integrals depending on parameters.
Introduce special functions 
\begin{equation*}
\sigma_{\la}(t,\tau)=
\frac{\thi_{1}(\la-t,\tau)\thi_{1}'(0,\tau)}
{\thi_{1}(\la,\tau)\thi_{1}(t,\tau)},
\qquad
E(t,\tau)=\frac{\thi_{1}(t,\tau)}{\thi_{1}'(0,\tau)}. 
\end{equation*}
They have the properties
\begin{gather}\label{tp1}
E(t+1,\tau)=-E(t,\tau),\qquad E(t+\tau,\tau)=-e^{-\pi i\tau-2\pi i t}E(t,\tau),\\
\label{tp2}
\sigma_{\la}(t+1,\tau)=\sigma_{\la}(t,\tau),
\qquad
\sigma_{\la}(t+\tau,\tau)=e^{2\pi i\la}\sigma_{\la}(t,\tau).
\end{gather}
The modular transformation properties of the functions $\sigma_{\la}$
and $E$ are calculated using the modular transformation properties of
the function $\thi_{1}$,
\begin{gather}\label{thetatrans}
\thi_{1}(t,\tau+1)=e^{\frac{\pi i}{4}}\thi_{1}(t,\tau),\qquad
\thi_{1}\left(\frac{t}{\tau},-\frac{1}{\tau}\right)=\frac{\sqrt{-i\tau}}{i}e^{\pi i\frac{t^2}{\tau}}\thi_{1}(t,\tau),
\end{gather}
where $\sqrt{-i\tau}$ is to be interpreted by the convention $|\arg(-i\tau)|<\pi/2$. 
We have
\begin{gather*}
\sigma_{\la}(t,\tau+1)=\sigma_{\la}(t,\tau),\qquad
\sigma_{\frac{\la}{\tau}}\left(\frac{t}{\tau},-\frac{1}{\tau}\right)=\tau
e^{-\frac{2\pi it\la}{\tau}}\sigma_{\la}(t,\tau),\\
E(t,\tau+1)=E(t,\tau),\qquad
E\left({t\over\tau},-{1\over\tau}\right)=\tau^{-1}e^{\pi i
  {t^{2}\over\tau}}E(t,\tau). 
\end{gather*}
Let $\Phi_{p,\ka}$ be the multi-valued function
\begin{equation*}
\Phi_{p,\ka}=\Phi_{\ka}(t_{1},\dots,t_{p},\tau)=
\prod_{j=1}^{p}E(t_{j},\tau)^{-{2p\over\ka}}\, 
\prod_{1\leq i<j\leq p}E(t_i - t_j ,\tau)^{2\over\ka}.
\end{equation*}
For $\ka\geq2$ and $n\in\mathbb{Z}$, let $\theta_{\ka,n}$ be the
theta function of level $\ka$ 
\begin{equation*}
\theta_{\ka,n}(t,\tau)=
\sum_{j\in\mathbb{Z}}e^{2\pi i\ka(j+\frac{n}{2\ka})^{2}\tau+2\pi
  i\ka(j+\frac{n}{2\ka})t}.
\end{equation*}
We have
\begin{gather}\label{reflp}
\theta_{\ka,n+2\ka}(t,\tau)=\theta_{\ka,n}(t,\tau),\qquad
\theta_{\ka,n}(-t,\tau)=\theta_{\ka,-n}(t,\tau),\\\label{tp3}
\theta_{\ka,n}\left(t + {2\over\ka},\tau\right)=q^{2n}\,\theta_{\ka,n}(t,\tau),\qquad
\theta_{\ka,n}\left(t + {2\tau\over\ka},\tau\right)=e^{-2\pi it-{2\pi
    i\tau\over\ka}}\,\theta_{\ka,n+2}(t,\tau).
\end{gather}
The modular transformation properties of $\theta_{\ka,n}$ are given by
\begin{gather*}
\theta_{\ka,n}(t,\tau+1)=q^{\frac{n^{2}}{2}}\theta_{\ka,n}(t,\tau),\\
\theta_{\ka,n}\left({t\over\tau},-{1\over\tau}\right)=\sqrt{-{i\tau\over2\ka}}e^{\pi
    i\ka{t^{2}\over2\tau}}\,\sum_{m=0}^{2\ka-1}q^{-mn}\,\theta_{\ka,m}(t,\tau).
\end{gather*}
Let $\Delta_{k}\subset\mathbb{R}^{k}\subset\mathbb{C}^{k}$ be the simplex 
\begin{equation*}
\Delta_{k}=\{(t_{1},t_{2},\dots ,t_{k})\in\mathbb{R}^{k}\subset\mathbb{C}^{k}\,|\,0\leq t_{k}\leq
    t_{k-1}\leq\cdots\leq t_{1}\leq1\}.
\end{equation*}
Let
$\widetilde{\Delta}_{k}$ be the image of ${\Delta_{k}}$ under the map $(t_{1},\dots,t_{k})\mapsto(\tau
t_{1},\dots,\tau t_{k})$. 
For $0\leq k\leq p$,
define 
\begin{equation*}
J_{\ka,n}^{[k]}(\la,\tau)=\int\Phi_{\ka}(t_{1},\dots,t_{p},\tau)\theta_{\ka,n}\left(\la+{2\over\ka}\sum_{j=1}^{p}t_{j},\tau\right)\,\prod_{j=1}^{p}\sigma_{\la}(t_{j},\tau)dt_{j}\,,
\end{equation*}
where the integral is over 
\begin{equation*}
\{(t_{1},t_{2},\dots,t_{p})\in\mathbb{C}^{p}\,|\,(t_{1},\dots,t_{k})\in\Delta_{k},\,(t_{k+1},\dots,t_{p})\in\widetilde{\Delta}_{p-k}\}.
\end{equation*}
The integral is a meromorphic function of the exponents of
$\Phi_{p,\ka}$.  It is well-defined when all of the exponents in the function
$\Phi_{p,\ka}$ have positive real parts.  We are interested in the
case when the exponents in the function $\Phi_{p,\ka}$ are negative
real numbers.  In this case, the integral is understood as an analytic
continuation from the region where the exponents have positive real parts.      
The branch of $\Phi_{p,\ka}$ is determined by fixing the arguments of all
factors of $\Phi_{p,\ka}$ for the case $\tau=i$ and deforming
continuously for arbitrary values of $\tau$.  For $\tau=i$, we fix $\arg E(t_{j},\tau)=0$
for $j=1,\dots,k$, $\arg E(t_{j},\tau)=\frac{\pi}{2}$ for
$j=k+1,\dots,p$, and $\arg E(t_{i}-t_{j},\tau)\in(-\pi,\pi)$ for $1\leq
i<j\leq p$.  
Introduce
\begin{equation*}
u_{n}^{[k]}(\la,\tau)=u^{[k]}_{\ka,n}(\la,\tau)=J_{\ka,n}^{[k]}(\la,\tau)+(-1)^{p+1}
J_{\ka,n}^{[k]}(-\la,\tau).
\end{equation*}
\begin{lemma}\label{l1}
The integrals $u_{n}^{[k]}$ have the properties
\begin{equation*}
u_{n}^{[k]}=u_{n+2\ka}^{[k]},\qquad
u^{[k]}_{n}=-q^{2k(n+p-k)}u^{[k]}_{-n-2(p-k)}.
\end{equation*}
\end{lemma}
\begin{proof}
The first equation is an immediate consequence of the $2\ka$
periodicity of the functions $\theta_{\ka,n}$.  To derive the second
equation, we consider the integrals $u_{n}^{[k]}$ after the change of variables
$t_{j}\mapsto-t_{j}+1$ for $1\leq j\leq k$, $t_{j}\mapsto-t_{j}+\tau$
for $k+1\leq j\leq p$.  The result follows from the formulas
\Ref{tp1}, \Ref{tp2}, \Ref{reflp}, and \Ref{tp3}. 
\end{proof}
\begin{thm}\cite{FV1}
For $0\leq k\leq p$ and any $n$, the integrals $u_{\ka,n}^{[k]}(\la,\tau)$ are solutions of the KZB-heat
equation having the properties (i)-(iv).    
\end{thm}
\begin{thm}  The set
\begin{equation*} 
 \{u_{n}^{[p]}(\la,\tau)\,\vert\, p+1\leq n\leq \ka -p-1\}
\end{equation*} 
is a basis for the space of conformal blocks.  
\end{thm}
\begin{proof}
We prove that for $n\in\{p+1,\dots,\ka-p-1\}$, the integrals $u_{n}^{[p]}$
are linearly independent over $\mathbb{C}$. 
 In Theorem \ref{zerothm}, we show that
for all other values of $n$ in the interval $0\leq n\leq \ka$, the
integrals $u_{n}^{[p]}$ are identically zero.  In the limit as $\tau\to i\infty$, the leading term of
$u_{n}^{[p]}$ is of the form
\begin{equation*} 
\frac{A_{p}e^{\frac{\pi i n^{2}}{2\ka}\tau}}{(\sin(\pi\la))^{p}}
\left(\left(\prod_{j=1}^{p}(q^{2(n+j)}-1)\right)
B_{p}\left(\frac{n+1}{\ka},-\frac{2p}{\ka},\frac{1}{\ka}\right)
\sin(\pi\la(n+p))+C_{n}(e^{\pi i\la},e^{-\pi i\la})\right),  
\end{equation*}
where $C_{n}(e^{\pi i\la},e^{-\pi i\la})$ is a Laurent
polynomial of degree $< n+p$ in $e^{\pi i\la}$.  Here $A_{p}$ is a non-zero constant depending only on $p$,
and $B_{p}(\alpha,\beta,\gamma)$ is the Selberg integral,
\begin{equation*}
B_{p}(\alpha,\beta,\gamma)=
{1\over p!}
\prod_{j=0}^{p-1}
\frac{\Gm(1+\gamma+j\gamma)
\Gm(\alpha+j\gamma)
\Gm(\beta+j\gamma)}
{\Gm(1+\gamma)
\Gm(\alpha+\beta+(p+j-1)\gamma)}.
\end{equation*}
It
follows that the $u_{n}^{[p]}$ are linearly independent provided that the
coefficients $\prod_{j=1}^{p}(q^{2(n+j)}-1)$ and
$B_{p}\left((n+1)/\ka,-2p/\ka,1/\ka\right)$ are nonzero functions for each $n$ in the range
$p+1\leq n\leq\ka-p-1$.  It is straightforward to check that the product
$\prod_{j=1}^{p}(q^{2(n+j)}-1)$ is never zero for $n$ in this
interval.  The Selberg integral
$B_{p}\left((n+1)/\ka,-2p/\ka, 1/\ka\right)$ is
also nonzero for $n$ in this interval since 
$1+1/\ka\notin\mathbb{Z} _{\leq 0}$, and
$(n-p+j)/\ka\notin\mathbb{Z} _{\leq0}$ for any $j$ satisfying $0\leq j\leq p-1$.     
\end{proof}  
We use the modular transformation properties of the functions $E$,
$\sigma_{\lambda}$ and $\theta_{\ka,n}$ to
obtain the following two lemmas.
\begin{lemma}\label{ttrans}
For any $n$, the action of the operator $T$ on $u_{n}^{[p]}$ is given
by the formula
\begin{equation*}
Tu_{n}^{[p]}(\la,\tau)=q^{\frac{n^{2}}{2}}u_{n}^{[p]}(\la,\tau).\qed
\end{equation*} 
\end{lemma}
\begin{lemma}\label{strans} For any $n$, we have
\begin{equation*} 
Su_{n}^{[p]}(\la,\tau)=\frac{e^{-\frac{\pi i}{4}}}{\sqrt{2\ka}}\sum_{m=0}^{2\ka-1}q^{-mn}u^{[0]}_{m}(\la,\tau). \qed
\end{equation*}
\end{lemma}
In Lemma \ref{strans}, we have written $Su_{n}^{[p]}$ as a linear
combination of $u_{m}^{[0]}$.  
Our goal is to express $Su_{n}^{[p]}$ in terms of the basis $ 
\{u_{m}^{[p]}(\la,\tau)\,\vert\, p+1\leq m\leq \ka -p-1\}$.   
This is accomplished using the Stokes theorem.  Repeated applications
of the Stokes theorem give us a recursive procedure for expressing 
the integrals $u_{m}^{[0]}$
as linear combinations of the integrals $u_{m}^{[p]}$.
 \begin{lemma}\label{st}For $0\leq k\leq p-1$ and any $n$, we have
\begin{equation*}
[p-k](q^{n+p-k}-q^{-n-p+k})u^{[k]}_{ n}=q^{-n-k-1}[k+1]
\left(q^{-2(k+1)}u^{[k+1]}_{n+2}-u^{[k+1]}_{n}\right).
\end{equation*}
\end{lemma}
\begin{proof}
Recall that the
$u_{n}^{[k]}$ are defined as analytic continuations of integrals where
all of the exponents in the function $\Phi_{p,\ka}$ have positive real
parts.  So the identity in the lemma relates objects which are
understood as analytic continuations.  
To prove the identity, we begin by considering simpler objects.
Namely, we consider the case when all of the exponents in 
$\Phi_{p,\ka}$ have positive real parts.  In this case, we can apply
the Stokes theorem which gives an identity for integrals with positive
exponents.  Then, we analytically continue all terms of the identity
to get the statement of the lemma.  More precisely, 
let $A\subset\mathbb{C}$ be the parallelogram with vertices at
$0,\,1,\,\tau,\,1+\tau$.  Consider the $(p+1)$-dimensional cell
\begin{equation*}
B_{k}=\{(t_{1},t_{2},\dots,t_{p})\in\mathbb{C}^{p}\,|\,(t_{1},\dots,t_{k})\in\Delta_{k},\,t_{k+1}\in
A,\,(t_{k+2},\dots,t_{p})\in\widetilde{\Delta}_{p-k-1}\}.
\end{equation*}
Applying the Stokes theorem to $B_{k}$ gives
\begin{multline}\label{stidentity}
\int_{\delta
  B_{k}}\Phi_{\ka}(t_{1},\dots,t_{p},\tau)\theta_{\ka,n}\left(\la+{2\over\ka}\sum_{j=1}^{p}t_{j},\tau\right)\prod_{j=1}^{p}\sigma_{\la}(t_{j},\tau)dt_{j}\\
+(-1)^{p+1}\int_{\delta
  B_{k}}\Phi_{\ka}(t_{1},\dots,t_{p},\tau)\theta_{\ka,n}\left(-\la+{2\over\ka}\sum_{j=1}^{p}t_{j},\tau\right)\prod_{j=1}^{p}\sigma_{-\la}(t_{j},\tau)dt_{j}=0.
\end{multline}  
The boundary of $B_{k}$ consists of $2p+2$
components of dimension $p$.  When all of the exponents in
the function $\Phi_{p,\ka}$ have positive real parts, the restrictions
of the integrand to all but four of the boundary
components are zero.  Those four components are
\begin{align*}
&\gamma^{1}_{k}=\{(t_{1},t_{2},\dots,t_{p})\in\mathbb{C}^{p}\,|\,(t_{1},\dots,t_{k})\in\Delta_{k},\,t_{k+1}\in
[0,1],\,(t_{k+2},\dots,t_{p})\in\widetilde{\Delta}_{p-k-1}\},\\
&\gamma^{2}_{k}=\{(t_{1},t_{2},\dots,t_{p})\in\mathbb{C}^{p}\,|\,(t_{1},\dots,t_{k})\in\Delta_{k},\,t_{k+1}\in
[1,1+\tau],\,(t_{k+2},\dots,t_{p})\in\widetilde{\Delta}_{p-k-1}\},\\
&\gamma^{3}_{k}=\{(t_{1},t_{2},\dots,t_{p})\in\mathbb{C}^{p}\,|\,(t_{1},\dots,t_{k})\in\Delta_{k},\,t_{k+1}\in
[\tau,1+\tau],\,(t_{k+2},\dots,t_{p})\in\widetilde{\Delta}_{p-k-1}\},\\
&\gamma^{4}_{k}=\{(t_{1},t_{2},\dots,t_{p})\in\mathbb{C}^{p}\,|\,(t_{1},\dots,t_{k})\in\Delta_{k},\,t_{k+1}\in
[0,\tau],\,(t_{k+2},\dots,t_{p})\in\widetilde{\Delta}_{p-k-1}\},
\end{align*}
where $[a,b]$ denotes the straight line segment connecting $a$ and
$b$.  Thus, the Stokes theorem applied to $B_{k}$ gives us that the sum
of the integrals over $\gamma_{1}^{k}$, $\gamma_{2}^{k}$,
$\gamma_{3}^{k}$, and $\gamma_{4}^{k}$ is equal to zero.

We obtain the result of the lemma as follows.  Let $B_{k}$ be defined
as above.  Take the analytic continuation of all terms in the identity
\Ref{stidentity}.  Let $\gamma^{*}_{k}$ be any component of
$\delta B_{k}\backslash\cup_{i=1}^{4}\gamma^{i}_{k}$.  In the case
when the exponents in $\Phi_{p,\ka}$ had positive real parts, we had
that the integral over $\gamma^{*}_{k}$ was zero.  Thus we must have
that the analytic continuation of the integral over $\gamma^{*}_{k}$
is also zero.  So we have that the sum of the integrals (understood
as analytic continuations) over $\gamma_{1}^{k}$, $\gamma_{2}^{k}$,
$\gamma_{3}^{k}$, and $\gamma_{4}^{k}$ is equal to zero.     
The integrals over the boundary components $\gamma_{1}^{k}$ and
$\gamma_{4}^{k}$ are $c^{k}_{1}u_{n}^{[k+1]}$ and
$c^{k}_{4}u_{n}^{[k]}$, respectively, where $c_{1}^{k},c_{4}^{k}\in\mathbb{C}$.  To recognize the integral over
$\gamma_{2}^{k}$ as $c_{2}^{k}u_{n}^{[k]}$, where $c_{2}^{k}\in\mathbb{C}$, we make the change of
variables $t_{k+1}\to t_{k+1}+1$.  To recognize the
integral over
$\gamma_{3}^{k}$ as $c_{3}^{k}u_{n+2}^{[k+1]}$, where $c_{3}^{k}\in\mathbb{C}$, we make the change of
variables $t_{k+1}\to t_{k+1}+\tau$.  
We calculate the constants
$c_{j}^{k}$ using the formulas \Ref{tp1}, \Ref{tp2}, and \Ref{tp3}.
\end{proof} 
\begin{thm} \label{zerothm} For $0\leq k\leq p$ and
  $n\in\{-p,-p+1,\dots,-p+2k\}\cup\{\ka-p,\ka-p+1,\dots,\ka-p+2k\}$,
  we have 
\begin{equation*}
u_{n}^{[k]}=0.
\end{equation*}
\end{thm}
\begin{proof}
We prove the statement by induction on $k$.  Setting $k=0$ in Lemma \ref{l1} gives
us the identity 
\begin{equation*}
u_{n}^{[0]}=-u_{-n-2p}^{[0]}.
\end{equation*}
Hence, using the $2\ka$ periodicity of the $u_{n}^{[k]}$, for
$n\equiv-p\mod\ka$, we have 
\begin{equation*}
u_{n}^{[0]}=0.
\end{equation*} 
Now assume the theorem is true for some
$k$ in the interval from $0$ to $p-1$.  We must show that $u_{n}^{[k+1]}=0$ for
$n\in\{-p,-p+1,\dots,-p+2k+2\}\cup\{\ka-p,\ka-p+1,\dots,\ka-p+2k+2\}$.
We use the reflection identity of Lemma \ref{l1},
\begin{gather}\notag 
u^{[k+1]}_{n}=-q^{2(k+1)(n+p-k-1)}u^{[k+1]}_{-n-2(p-k-1)},
\intertext{and the relation of Lemma \ref{st},}
\label{st2}
[p-k](q^{n+p-k}-q^{-n-p+k})u^{[k]}_{n}=q^{-n-k-1}[k+1]
\left(q^{-2(k+1)}u^{[k+1]}_{n+2}-u^{[k+1]}_{n}\right).
\end{gather}
By the induction hypothesis, the left hand side of
equation \Ref{st2} is zero for
$n\in\{-p,-p+1,\dots,-p+2k\}\cup\{\ka-p,\ka-p+1,\dots,\ka-p+2k\}$.
Thus, for
$n\in\{-p,-p+1,\dots,-p+2k\}\cup\{\ka-p,\ka-p+1,\dots,\ka-p+2k\}$, we
have the recursion relations
\begin{equation}\label{rr}
u_{n}^{[k+1]}=q^{-2(k+1)}u_{n+2}^{[k+1]}.
\end{equation}
We use the identities
\begin{equation*}
u^{[k+1]}_{-p+k+1}=u^{[k+1]}_{\ka-p+k+1}=0, \quad u^{[k+1]}_{-p+k}=u^{[k+1]}_{\ka-p+k}=0,
\end{equation*}
to obtain, recursively on the subscripts, the result.  The first identity is an
immediate consequence of the reflection identity, and the second
identity follows from comparing the recursion relations \Ref{rr} for $n=-p+k$ and
$n=\ka-p+k$, respectively, with the reflection identities 
\begin{equation*}
u^{[k+1]}_{-p+k}=-q^{-2(k+1)}u^{[k+1]}_{-p+k+2}, \qquad u^{[k+1]}_{\ka-p+k}=-q^{-2(k+1)}u^{[p]}_{\ka-p+k+2}.
\end{equation*}
\end{proof}
\begin{proposition}\label{5c0}  For $0\leq k\leq p$ and any $n$,
\begin{equation*}
Su_{n}^{[p]}(\la,\tau)=\frac{e^{-\frac{\pi i}{4}}}{\sqrt{2\ka}}\left(\sum_{m=-p+2k+1}^{\ka-p-1}f^{(k)}_{m,n}u_{m}^{[k]}(\lambda,\tau)+\sum_{m=\ka-p+2k+1}^{2\ka-p-1}f^{(k)}_{m,n}u_{m}^{[k]}(\lambda,\tau)\right),
\end{equation*}
where
\begin{equation*}
f_{m,n}^{(k)}=\frac{q^{-m(n+k)-\frac{k(k+1)}{2}}}{(q^{-1}-q)^{k}}\begin{bmatrix}p\\k\end{bmatrix}^{-1}
\sum_{j=0}^{k}\begin{bmatrix}k\\j\end{bmatrix}\frac{q^{2jn}}{(-m-p+k+1,q)_{j}\,(m+p-k+1,q)_{k-j}}.
\end{equation*}
\end{proposition}
\begin{proof}
The proof is by induction on $k$.  The result from Lemma
\ref{strans} (where we have used the periodicity of the $u_{n}^{[k]}$
to shift the interval of summation) combined with the identity
$u_{-p}^{[0]}=u_{\ka-p}^{[0]}=0$ of Theorem \ref{zerothm} give us  
\begin{equation*} 
e^{\frac{\pi i}{4}}\sqrt{2\ka}\,Su_{n}^{[p]}(\la,\tau)=S'u_{n}^{[p]}=\sum_{m=-p+1}^{\ka-p-1}q^{-mn}u^{[0]}_{m}(\la,\tau)+
\sum_{m=\ka-p+1}^{2\ka-p-1}q^{-mn}u^{[0]}_{m}(\la,\tau).
\end{equation*}
But this is exactly the result for $k=0$ since 
\begin{equation*}
f^{(0)}_{m,n}=q^{-mn}.
\end{equation*}
Assuming the result for some $k$ satisfying $0\leq k\leq p-1$, we have
\begin{equation}\label{ih}
S'u_{n}^{[p]}(\la,\tau)=\sum_{m=-p+2k+1}^{\ka-p-1}f^{(k)}_{m,n}u_{m}^{[k]}(\lambda,\tau)+\sum_{m=\ka-p+2k+1}^{2\ka-p-1}f^{(k)}_{m,n}u_{m}^{[k]}(\lambda,\tau).
\end{equation}
Rewriting \Ref{ih} using the identity of Lemma \ref{st},
\begin{equation*}
u^{[k]}_{m}=\frac{q^{-m-k-1}}{(q^{m+p-k}-q^{-m-p+k})}\frac{[k+1]}{[p-k]}
\left(q^{-2(k+1)}u^{[k+1]}_{m+2}-u^{[k+1]}_{m}\right),
\end{equation*}
which makes sense for all $m$ in the interval of summation since
those values of $m$ never give zero denominators, we obtain
\begin{multline}\label{ih2}
S'u_{n}^{[p]}(\la,\tau)=
\frac{q^{-k-1}}{(q^{-1}-q)}\frac{[k+1]}{[p-k]}\,\times\\
\biggl(\,
\sum_{m=-p+2k+1}^{\ka-p-1}
\frac{q^{-m}f^{(k)}_{m,n}}{[m+p-k]}u_{m}^{[k+1]}(\lambda,\tau)
-
\sum_{m=-p+2k+3}^{\ka-p+1}
\frac{q^{-m-2k}f^{(k)}_{m-2,n}}{[m+p-k-2]}u_{m}^{[k+1]}(\lambda,\tau)\\
+
\sum_{m=\ka-p+2k+1}^{2\ka-p-1}
\frac{q^{-m}f^{(k)}_{m,n}}{[m+p-k]}u_{m}^{[k+1]}(\lambda,\tau)
-
\sum_{m=\ka-p+2k+3}^{2\ka-p+1}
\frac{q^{-m-2k}f^{(k)}_{m-2,n}}{[m+p-k-2]}u_{m}^{[k+1]}(\lambda,\tau)
\biggl)
.
\end{multline}
By Theorem \ref{zerothm}, the integrals $u_{m}^{[k+1]}$ are identically
zero for $m\in\{-p+2k+1,-p+2k+2,\ka-p,\ka-p+1,\ka-p+2k+1,\ka-p+2k+2,2\ka-p,2\ka-p+1\}$.
Thus, we may combine the first two terms on the right hand side of
\Ref{ih2} into one sum over $m$ in the range $-p+2k+3\leq m \leq \ka-p-1$.
Similarly, the last two terms give us a sum over $m$ in the range $\ka-p+2k+3
\leq m \leq 2\ka-p-1$.  In each of these sums, the coefficient of
$u_{m}^{[k+1]}$ is given by the expression
\begin{multline}\label{coeff}
\frac{q^{-m(n+k+1)-\frac{(k+1)(k+2)}{2}}}{(q^{-1}-q)^{k+1}}
\begin{bmatrix}p\\k+1\end{bmatrix}^{-1}\biggl(
\sum_{j=0}^{k}\begin{bmatrix}k\\j\end{bmatrix}
\frac{q^{2jn}}{(-m-p+k+1,q)_{j}\,(m+p-k,q)_{k-j+1}}\\
+
\sum_{j=1}^{k+1}\begin{bmatrix}k\\j-1\end{bmatrix}
\frac{q^{2jn}}{(-m-p+k+2,q)_{j}\,(m+p-k-1,q)_{k-j+1}}\biggl).
\end{multline}
To complete the proof of the theorem, we must show that the expression \Ref{coeff}
is equal to the expression
\begin{equation*}
\frac{q^{-m(n+k+1)-\frac{(k+1)(k+2)}{2}}}{(q^{-1}-q)^{k+1}}\begin{bmatrix}p\\k+1\end{bmatrix}^{-1}
\sum_{j=0}^{k+1}\begin{bmatrix}k+1\\j\end{bmatrix}\frac{q^{2jn}}{(-m-p+k+2,q)_{j}\,(m+p-k,q)_{k-j+1}}.
\end{equation*}
This is proved by direct calculation.
\end{proof}
\begin{corollary}\label{5c1}
For $0\leq k\leq p-1$, $m\in\{-p+2k+1,\dots,\ka-p-1\}\cup\{\ka-p+2k+1,\dots,2\ka-p-1\}$, and
any $n$, we have
\begin{equation*}
f_{m,n}^{(k+1)}=\frac{q^{-m-k-1}}{q-q^{-1}}\,\frac{[k+1]}{[p-k]}\left(\frac{q^{-2k}f^{(k)}_{ m-2,n}}{[m-2+p-k]}-\frac{f^{(k)}_{ m,n}}{[m+p-k]}\right).
\end{equation*}
\end{corollary}
\begin{lemma}\label{l2}
The functions $f^{(k)}_{m,n}$ have the property
\begin{equation*}
f^{(k)}_{m,n}=q^{-2k(m+p-k)+2pn}f^{(k)}_{-m-2p+2k,-n}.\qed
\end{equation*}
\end{lemma}
We rewrite the sum in Proposition \ref{5c0} using Lemmas \ref{l1} and \ref{l2}.
\begin{corollary}
For $0\leq k\leq p$ and any $n$,
\begin{equation*}
Su^{[p]}_{
  n}(\la,\tau)=\frac{e^{-\frac{\pi i}{4}}}{\sqrt{2\ka}}\sum_{m=-p+2k+1}^{\ka-p-1}
\left(f_{m,n}^{(k)}-q^{2pn}f_{m,-n}^{(k)}\right)u_{
  m}^{[k]}(\la,\tau).
\end{equation*}
\end{corollary}
\section{Macdonald Polynomials and the Shift Operator}
The Macdonald polynomials of type $A_{1}$ are $x$-even polynomials
in terms of $q^{mx}$, where $m\in\mathbb{Z}$.  They depend on two
parameters $k$ and $n$, where $k$ and $n$ are non-negative
integers.  They are
defined by the conditions:
\begin{enumerate}  
\item 
$P_{n}^{(k)}(x)=q^{nx}+q^{-nx}+$ lower order terms, except for $P_{0}^{(k)}(x)=1$,     
\item
$\langle P_{m}^{(k)},P_{n}^{(k)}\rangle=0$ for $m\neq n$, where
\begin{equation*}
\langle f,g\rangle= \,\frac{1}{2}\text{Const Term}\,\left(fg\prod_{j=0}^{k-1}(1-q^{2(j+x)})(1-q^{2(j-x)})\right).
\end{equation*}
\end{enumerate}
\begin{example}  For $n>0$, we have
\begin{equation*}
P_{n}^{(0)}(x)=q^{nx}+q^{-nx}.
\end{equation*}
\end{example}

The shift operator $D$ is an operator acting on functions $f(x)$ by
\begin{equation*} 
Df(x)=\frac{f(x-1)-f(x+1)}{q^{x}-q^{-x}}.
\end{equation*}
The name is ``shift operator'' because its action on the basic
($q$-difference) hypergeometric functions results in a shift of the
parameters \cite{AW,Ch}.     
\begin{thm}\cite{AI}\label{so}  
For $n\geq 1$ and $k\geq 0$, we have
\begin{equation*}
DP_{n}^{(k)}(x)=(q^{-n}-q^{n})\,P_{n-1}^{(k+1)}(x).
\end{equation*}
\end{thm}

\begin{remark}
It follows from Theorem \ref{so} that all of the Macdonald polynomials
can be calculated from the polynomials $P^{(0)}_{n}$ using the shift
operator provided that $q^{-n}-q^{n}\neq0$.
\end{remark}
\begin{example} For $n>0$, we have
\begin{equation*}
P^{(1)}_{n-1}(x)=\frac{q^{nx}-q^{-nx}}{q^x-q^{-x}}.
\end{equation*}
\end{example} 
\section{Identification of the $f_{m,n}^{(k)}$ with values of the
  Macdonald polynomials}
\begin{thm}\label{mp=f}  For $0\leq k\leq p$, $k+1\leq n \leq \ka$, and
  $m\in\{-p+2k+1,\dots,\ka-p-1\}$, we have
\begin{multline*}
f_{m,n}^{(k)}-q^{2pn}f_{m,-n}^{(k)}= q^{pn-km-\frac{k(k+1)}{2}}\begin{bmatrix}p\\k\end{bmatrix}^{-1}(q^{-m-p+k}-q^{m+p-k})\;\times\\
\left(\prod_{j=1}^{k}(q^{-n+j}-q^{n-j})\right)P^{(k+1)}_{n-k-1}(m+p-k).
\end{multline*}
\end{thm}
\begin{proof}
We prove the statement by induction on $k$.  The case $k=0$ follows directly
from the formulas
\begin{equation*}
f_{m,n}^{(0)}=q^{-mn},
\qquad
P^{(1)}_{n-1}(m+p)=\frac{q^{n(m+p)}-q^{-n(m+p)}}{q^{m+p}-q^{-m-p}}.  
\end{equation*}

We assume that the theorem is true for some $k$ in the range $0\leq k\leq p-1$.
To prove the theorem, we must verify the identity 
\begin{multline}\label{i1}
f_{m,n}^{(k+1)}-q^{2pn}f_{m,-n}^{(k+1)}=
q^{pn-(k+1)m-\frac{(k+1)(k+2)}{2}}
\begin{bmatrix}p\\k+1\end{bmatrix}^{-1}
(q^{-m-p+k+1}-q^{m+p-k-1})\,\times\\
\left(\prod_{j=1}^{k+1}(q^{-n+j}-q^{n-j})\right)P^{(k+2)}_{n-k-2}(m+p-k-1)
\end{multline} 
for $k+2\leq n\leq \ka$ and 
$m\in\{-p+2k+3,\dots,\ka-p-1\}$.
By Corollary \ref{5c1}, the left hand side of the equation \Ref{i1} equals
\begin{equation*}
\frac{q^{-m-k-1}}{(q-q^{-1})}\frac{[k+1]}{[p-k]}\left(q^{-2k}\,\frac{f_{m-2,n}^{(k)}-q^{2pn}f_{m-2,-n}^{(k)}}{[m+p-k-2]}-\frac{f_{m,n}^{(k)}-q^{2pn}f_{m,-n}^{(k)}}{[m+p-k]}\right).
\end{equation*}
By Theorem \ref{so}, the right hand side of the equation \Ref{i1} equals
\begin{multline*}
 q^{pn-(k+1)m-\frac{(k+1)(k+2)}{2}}\begin{bmatrix}p\\k+1\end{bmatrix}^{-1}
\left(\prod_{j=1}^{k}(q^{-n+j}-q^{n-j})\right)\,\times\\
\left(P^{(k+1)}_{n-k-1}(m+p-k)-P^{(k+1)}_{n-k-1}(m+p-k-2)\right).
\end{multline*}
Thus, the equation \Ref{i1} is equivalent to the equation
\begin{multline*}
q^{-2k}\,\frac{f_{m-2,n}^{(k)}-q^{2pn}f_{m-2,-n}^{(k)}}{q^{m+p-k-2}-q^{-m-p+k+2}}+\frac{f_{m,n}^{(k)}-q^{2pn}f_{m,-n}^{(k)}}{q^{-m-p+k}-q^{m+p-k}}
=
q^{pn-km-\frac{k(k+1)}{2}}\begin{bmatrix}p\\k\end{bmatrix}^{-1}\,\times\\
\left(\prod_{j=1}^{k}(q^{-n+j}-q^{n-j})\right)\left(P^{(k+1)}_{n-k-1}(m+p-k)-P^{(k+1)}_{n-k-1}(m+p-k-2)\right),
\end{multline*}
which follows from the induction hypothesis.
\end{proof}
\begin{corollary}\label{smf}
For $0\leq k\leq p$ and  $k+1\leq n \leq \ka$, 
\begin{multline*}
Su_{ n}^{[p]}(\la,\tau)=\frac{e^{\frac{-\pi i}{4}}}{\sqrt{2\ka}}\sum_{m=-p+2k+1}^{\ka-p-1}q^{pn-km-\frac{k(k+1)}{2}}\begin{bmatrix}p\\k\end{bmatrix}^{-1}(q^{-m-p+k}-q^{m+p-k})\;\times\\
\left(\prod_{j=1}^{k}(q^{-n+j}-q^{n-j})\right)P^{(k+1)}_{n-k-1}(m+p-k)u_{ m}^{[k]}(\la,\tau).
\end{multline*}
\end{corollary}
 
The left hand side of the equation in Corollary \ref{smf} is equal to zero if
$n=k+1,...,p$ or $n=\kappa-p, \kappa -p + 1,...,\kappa - k - 1$.
This gives $2(p-k)$ relations for the $\kappa - 2k-1$ possibly nonzero functions
$u^{[k]}_m$, $m=k+1,...,p-k-1$. 
\begin{lemma}
The $2(p-k)$ relations between the functions $u^{[k]}_m$,
$m=k+1,...,p-k-1$, given by Corollary \ref{smf} are linearly
independent and thus generate all linear relations between those functions.
\end{lemma}
\begin{proof}
The coefficients of the above relations form $2(p-k)$ columns of the
$S$-matrix corresponding to $\kappa'=\kappa$ and $p'=k$. 
That $S$-matrix is a non-degenerate matrix whose square is a root of
unity. Thus any set of its columns is linearly independent.
\end{proof}
\section{Representation of the Modular Group}
Let $T=(t_{m,n})$ and $S=(s_{m,n})$ be
the matrices of the transformations $T$ and $S$, respectively, with
respect to the basis
\begin{equation*} 
 \{u_{\ka,n}^{[p]}(\la,\tau)\,\vert\, p+1\leq n\leq \ka -p-1\}.
\end{equation*}  
Here the matrices $T$ and $S$ are defined by
\begin{gather*}
Tu_{\ka,n}^{[p]}=\sum_{m=p+1}^{\ka-p-1}t_{m,n}u_{\ka,m}^{[p]},\\
Su_{\ka,n}^{[p]}=\sum_{m=p+1}^{\ka-p-1}s_{m,n}u_{\ka,m}^{[p]}.
\end{gather*}
\begin{thm}
For $p+1\leq m,n\leq \ka-p-1$, we have
\begin{gather*}
t_{m,n}=q^{\frac{n^2}{2}}\delta_{mn},\\
s_{m,n}=
\frac{e^{-\frac{\pi i}{4}}}{\sqrt{2\ka}}q^{p(n-m)-\frac{p(p+1)}{2}}(q^{-m}-q^{m})\left(\prod_{j=1}^{p}(q^{-n+j}-q^{n-j})\right)P^{(p+1)}_{n-p-1}(m),
\end{gather*}
where $\delta_{mn}=1$ for $m=n$ and $0$ otherwise.
\end{thm}
The theorem follows directly from Lemma \ref{ttrans} and Theorem \ref{smf}.
\begin{example} Let $\kappa= 2p+2$. In this case, the only element of the S-matrix
is
\begin{equation*}
s_{p+1,\,p+1} = \frac{(-i)^{p+1}}{\sqrt{p+1}}e^{-\pi i\frac{p+1}{4}}\prod_{j=1}^{p}(q^{j}+q^{-j}).
\end{equation*}
On the other hand, according to \cite{FSV},
\begin{equation*}
u^{[p]}_{2p+2,p+1} =a_{p}
\thi_{1}(\la,\tau)^{p+1},
\end{equation*}
where $a_{p}$ is a constant, 
and hence,
\begin{equation*}
Su^{[p]}_{2p+2,p+1} = a_{p}(-i)^{p+1}e^{-\pi i\frac{p+1}{4}}\theta (\lambda,\tau)^{p+1}.
\end{equation*}
This, in particular, implies 
\begin{equation*}
\prod_{j=1}^p(q^j+q^{-j})=\sqrt{p+1}.
\end{equation*}
This formula and its relations to the classical Gauss sums can be
found in \cite{Ch}.
\end{example}                                               
The projective representation of the modular
group in the space of conformal blocks described in \cite{K} 
is given by the matrices $\tilde{T}=(\tilde{t}_{m,n})$
and $\tilde{S}=(\tilde{s}_{m,n})$, $p+1\leq m,n\leq \ka-p-1$,
where
\begin{gather*}
\tilde{t}_{m,n}=e^{-\frac{\pi i}{4}}q^{\frac{n^2}{2}}\delta_{mn},\\
\tilde{s}_{m,n}=
\frac{i}{\sqrt{2\ka}}q^{-\frac{p(p+1)}{2}}\left(\prod_{j=0}^{p}(q^{-m+j}-q^{m-j})\right)P^{(p+1)}_{n-p-1}(m).
\end{gather*}
\begin{proposition}
Let $D$ be the diagonal matrix such that 
\begin{equation*}
d_{j}=q^{pj}\,\prod_{l=1}^{p}(q^{-j+l}-q^{j-l}),\qquad p+1\leq j\leq \ka-p-1.
\end{equation*}
Then
\begin{equation*}
T=e^{\frac{\pi i}{4}}D^{-1}\tilde{T}D,\qquad  S=e^{-\frac{3\pi i}{4}}D^{-1}\tilde{S}D.\qed
\end{equation*}
\end{proposition}

{\bf Remark.} For the standard definition of conformal blocks
$v(\lambda,\tau)=\vartheta_1(\lambda,\tau)^{-1}u(\lambda,\tau)$
(see the remark in Sect.~2), the $S$ and
$T$ transformations are defined as
\[
\hat Tv(\lambda,\tau)=v(\lambda,\tau+1),\qquad 
\hat Sv(\lambda,\tau)=e^{-\pi i (\kappa-2){\lambda^{2}\over2\tau}}\tau^{-
\frac{p(p+1)}{\kappa}}\,v\left({\lambda\over\tau},-{1\over\tau}\right).
\]
Using the transformation rules \Ref{thetatrans} of
$\vartheta_1$ we see that 
the corresponding matrices $\hat T$, $\hat S$ are related to $T$, $S$
by
\[
\hat T=e^{-\frac{\pi i}{4}}T,\qquad \hat S=e^{\frac{3\pi i}{4}}S.
\]
Thus $\hat T$, $\hat S$ are conjugated to $\tilde T$, $\tilde S$.

\section{Trace Functions for $U_{q}(sl_{2})$}  
Let $q=e^{\frac{\pi i}{\ka}}$.  The quantum group $U_{q}(sl_{2})$ has
generators $E,F,q^{ch}$, where $c\in\mathbb{C}$, with relations
\begin{gather*}
q^{ch}q^{c'h}=q^{(c+c')h},\,\quad
q^{ch}Eq^{-ch}=q^{2c}E,\,\quad q^{ch}Fq^{-ch}=q^{-2c}F,\,\quad
EF-FE=\frac{q^{h}-q^{-h}}{q-q^{-1}},
\intertext{and comultiplication defined by}
\Delta(E)=E\otimes q^{h}+1\otimes
E,\qquad\Delta(F)=F\otimes1+q^{-h}\otimes F,\qquad
\Delta(q^{ch})=q^{ch}\otimes q^{ch}.
\end{gather*}
Identify weights for
$U_{q}(sl_{2})$ with complex numbers as follows.  Say that a
vector $v$ in a $U_{q}(sl_{2})$-module has weight $\nu\in\mathbb{C}$ if
$q^{h}v=q^{\nu}v$.  Let $M_{\mu}$ be the Verma module over $U_{q}(sl_{2})$ with highest
weight $\mu$, and let $v_{\mu}$ be its highest weight vector.  Let $k$
be a non-negative integer such that $\ka\geq2k+2$.  Let $U$ be
the irreducible finite dimensional representation of $U_{q}(sl_{2})$
of weight $2k$.  Let $U[0]$ denote the zero weight subspace of $U$.
Let $u\in U[0]$.  For generic $\mu$, let
$\varPhi_{\mu}^{u}:M_{\mu}\to M_{\mu}\otimes U$ be the intertwining
operator defined by
\begin{equation*}  
\varPhi_{\mu}^{u}v_{\mu}=v_{\mu}\otimes u +\frac{Fv_{\mu}}{[-\mu]}\otimes
Eu+\cdots+\frac{F^{j}v_{\mu}}{[j]!(-\mu,q)_{j}}\otimes
  E^{j}u+\cdots.
\end{equation*}  
Introduce an $\End(U[0])$-valued function $\psi^{(k)}(q,\nu,\mu)$ defined by
\begin{equation*}
\psi^{(k)}(q,\nu,\mu)u=\tr|_{M_{\mu}}(\varPhi_{\mu}^{u}q^{\nu h}).
\end{equation*}
Since $U[0]$ is one-dimensional, this function is a scalar function.
\begin{thm}\cite{EV} The function $\psi^{(k)}(q,\nu,\mu)$ is given by the formula
\begin{equation*}
\psi^{(k)}(q,\nu,\mu)=q^{\nu\mu}\sum_{j=0}^{k}(-1)^{j}q^{\frac{j(j-3)}{2}}(q-q^{-1})^{-j-1}\frac{[k+j]!}{[j]![k-j]!}\frac{q^{-j\mu-(j-1)\nu}}{\prod_{l=0}^{j-1}[\mu-l]\prod_{l=0}^{j}[\nu-l]}.
\end{equation*}
\end{thm}
Introduce renormalized trace functions $\Psi^{(k)}(q,\nu,\mu)$
defined by
\begin{equation*}
\Psi^{(k)}(q,\nu,\mu)=\prod_{j=1}^{k}\left(\frac{q^{\mu+1-j}-q^{-\mu-1+j}}{q^{\nu+j}-q^{-\nu-j}}\right)\psi^{(k)}(q,\nu,\mu).   
\end{equation*}
The function $\Psi^{(k)}$ is a holomorphic function of $\mu$. 
\section{Identification of the $f_{m,n}^{(k)}$ with values
  of the
  trace functions $\Psi^{(k)}$ }
\begin{thm}\label{tf}
For $0\leq k\leq p$, $m\in\{-p+2k+1,\dots,\ka-p-1\}\cup\{\ka-p+2k+1,\dots,2\ka-p-1\}$, and
any $n$, 
\begin{equation*}
f_{m,n}^{(k)}=q^{pn-km-k(k+1)}\,(q^{m+p-k}-q^{-m-p+k})\begin{bmatrix}p\\k\end{bmatrix}^{-1}\Psi^{(k)}(q^{-1},-m-p+k,-n-1).
\end{equation*}
\end{thm}
\begin{proof}
The statement of the theorem is equivalent to the identity
\begin{multline}\label{f=psi}
(q-q^{-1})^{k}\sum_{j=0}^{k}\begin{bmatrix}k\\j\end{bmatrix}(-m-p+k+j+1,q)_{k-j}(m+p-j+1,q)_{j}q^{2jn}\\
=q^{-\frac{k(k+1)}{2}+kn}\sum_{j=0}^{k}q^{-\frac{j(j-1)}{2}-j(m+p-k+n)}(q-q^{-1})^{j}\frac{[k+j]!}{[j]![k-j]!}\,\times\\
\prod_{l=j+1}^{k}(q^{m+p-k+l}-q^{-m-p+k-l})(q^{n+l}-q^{-n-l}).
\end{multline}
Let $x=q^{n}$.  Using the
$q$-binomial theorem, we rewrite the identity \Ref{f=psi} as
\begin{multline}\label{f=psi2}
(q-q^{-1})^{k}\sum_{j=0}^{k}\begin{bmatrix}k\\j\end{bmatrix}(-m-p+k+j+1,q)_{k-j}(m+p-j+1,q)_{j}x^{2j-k}\\
=q^{-k(k+1)}\sum_{j=0}^{k}q^{-j(m+p-k-1)}(q-q^{-1})^{j}\frac{[k+j]!}{[j]![k-j]!}\,\times\\
\prod_{l=j+1}^{k}(q^{m+p-k+l}-q^{-m-p+k-l})\sum_{i=0}^{k-j}(-1)^{k-j-i}q^{i(k+j+1)}\begin{bmatrix}k-j\\i\end{bmatrix}x^{2i-k}.
\end{multline}
For a fixed value of $j$ in the range $0\leq j\leq k$, the coefficient
of $x^{2j-k}$ on the left hand side of \Ref{f=psi2} is
\begin{equation*}
(q-q^{-1})^{k}\begin{bmatrix}k\\j\end{bmatrix}(-m-p+k+j+1,q)_{k-j}(m+p-j+1,q)_{j}.
\end{equation*}
The coefficient
of $x^{2j-k}$ on the right hand side of \Ref{f=psi2} is
\begin{multline*}
q^{-(k+1)(k-j)}(q-q^{-1})^{k}\begin{bmatrix}k\\j\end{bmatrix}(m+p-j+1,q)_{j}(-m-p+j,q)_{k-j}
\,\times\\\sum_{i=0}^{k-j}q^{-i(m+p-k-j-1)}\frac{(k+1,q)_{i}(j-k,q)_{i}}{[i]!(m+p-k+1,q)_{i}}.
\end{multline*}
We observe that the sum 
\begin{equation*}
\sum_{i=0}^{k-j}q^{-i(m+p-k-j-1)}\frac{(k+1,q)_{i}(j-k,q)_{i}}{[i]!(m+p-k+1,q)_{i}}
\end{equation*}
is the hypergeometric series
\begin{equation}\label{hgs}
_{2}\phi_{1}\biggl(\begin{array}{c}q^{-2(k+1)},\,q^{-2(j-k)}\\q^{-2(m+p-k+1)}\end{array};\,q^{-2},\,q^{-2(m+p-k-j)}\biggl).
\end{equation}
The series \Ref{hgs} is equal to
\begin{equation*}
q^{(k-j)(k+1)}\frac{(m+p-2k,q)_{k+1}}{(m+p-k-j,q)_{k+1}}.
\end{equation*}
(For a proof of this formula, see \cite{GR}).
Hence the coefficients of $x^{2j-k}$ on both sides of \Ref{f=psi2} are
equal.  This proves the theorem.
\end{proof}
Under the identification in Theorem \ref{tf}, Corollary
  \ref{5c1} takes the following form.
\begin{corollary}\label{so2} For $0\leq k\leq p-1$, $m\in\{-p+2k+1,\dots,\ka-p-1\}\cup\{\ka-p+2k+1,\dots,2\ka-p-1\}$,
  and any $n$,
\begin{equation*}
q^{-k-1}\Psi^{(k+1)}(q^{-1},m,n+k)=
\frac{\Psi^{(k)}(q^{-1},m-1,n+k)-\Psi^{(k)}(q^{-1},m+1,n+k)}{q^m - q^{-m}}.
\end{equation*} 
\end{corollary} 
\begin{remark}  The right hand side of the equation in Corollary
  \ref{so2} is the shift operator applied to
  $\Psi^{(k)}(q^{-1},m,n+k)$.  Thus, for $m$ as above, $k\geq0$, and
  any $n$, we have 
\begin{equation*}
D\Psi^{(k)}(q^{-1},m,n)=q^{-k-1}\,\Psi^{(k+1)}(q^{-1},m,n).
\end{equation*}
\end{remark}

Together with Theorem \ref{mp=f}, the identification in Theorem \ref{tf} 
gives us a formula for the Macdonald polynomials evaluated at roots of unity in terms of
values of the renormalized trace functions. 
\begin{corollary}\label{es}
For $\ka\geq 2k+2$, $k+1\leq n\leq \ka$, and
$m\in\{-p+2k+1,\dots,\ka-p-1\}$, we have
\begin{multline*}
\Psi^{(k)}(q^{-1},-m-p+k,n-1)-\Psi^{(k)}(q^{-1},-m-p+k,-n-1)\\=
P^{(k+1)}_{n-k-1}(m+p-k)\prod_{j=1}^{k}(q^{-n+2j}-q^{n}).
\end{multline*}
\end{corollary}
\begin{remark}
For a generalization of Corollary \ref{es}, see \cite{ES}.
\end{remark}

Using the identification in Theorem \ref{tf}, we have the following expression for the $S$ transformation in terms of
the renormalized trace functions.
\begin{corollary}
For $0\leq k\leq p$ and any $n$, we have
\begin{multline*}
Su_{ n}^{[p]}(\la,\tau)=\frac{e^{-\frac{\pi i}{4}}}{\sqrt{2\ka}}\sum_{m=-p+2k+1}^{\ka-p-1}q^{pn-km-k(k+1)}\begin{bmatrix}p\\k\end{bmatrix}^{-1}(q^{-m-p+k}-q^{m+p-k})\;\times\\
\left(\Psi^{(k)}(q^{-1},-m-p+k,n-1)-\Psi^{(k)}(q^{-1},-m-p+k,-n-1)\right)u_{ m}^{[k]}(\la,\tau)\,.
\end{multline*}
\end{corollary}

\end{document}